\documentclass[10pt]{article}
\font\smallit=cmti10

\usepackage[export]{adjustbox}
\usepackage{forest}
\usepackage{amssymb,amsmath,amsthm,enumitem} 
\usepackage{subcaption,booktabs}
\usepackage{amstext}
\usepackage{amsfonts}
\usepackage{array}
\usepackage{hyperref}
\usepackage{url}
\usepackage{indentfirst}
\usepackage{tikz}
\usepackage{booktabs}
\usepackage{filecontents}
\usepackage[margin=1.5in]{geometry}    
\usepackage{graphicx}
\usepackage[utf8]{inputenc}
\usepackage{hyperref}
\usepackage{caption}
\usetikzlibrary{shapes.multipart}
\tikzset{block/.style={
        font=\sffamily,
        draw=black,
        thin,
        fill=white!50,
        rectangle split,
        rectangle split horizontal,
        rectangle split parts=#1,
        outer sep=0pt},
        }
\tikzset{every label/.style = {font=\footnotesize,  text=red,
                               inner sep=1pt }}

\usepackage{enumerate} 
\usepackage[justification=centering]{caption}
\usepackage{caption}
\usetikzlibrary{arrows.meta}

\makeatletter

\renewcommand\section{\@startsection {section}{1}{\z@}
{-30pt \@plus -1ex \@minus -.2ex}
{2.3ex \@plus.2ex}
{\normalfont\normalsize\bfseries\boldmath}}

\renewcommand\subsection{\@startsection{subsection}{2}{\z@}
{-3.25ex\@plus -1ex \@minus -.2ex}
{1.5ex \@plus .2ex}
{\normalfont\normalsize\bfseries\boldmath}}

\renewcommand{\@seccntformat}[1]{\csname the#1\endcsname. }

\makeatother

\newtheorem{theorem}{Theorem}

\newtheorem{lemma}{Lemma}

\newtheorem{proposition}{Proposition}

\newtheorem{conjecture}{Conjecture}
\newtheorem{observation}{Observation}

\theoremstyle{definition}

\begin{document}

\begin{center}
\uppercase{On partitioning the set \lowercase{$[n] = \{1, \dots, n\}$}
into subsets of size at most \lowercase{$m$}  \\
such that all sums are powers of \lowercase{$m$}}

\vskip 20pt
{\bf Vladimir Gurvich}\\
{\smallit National Research University Higher School of Economics (HSE), Moscow, Russia and\\ 
RUTCOR, Rutgers University,  Piscataway, NJ, United States}\\
{\tt vgurvich@hse.ru}, {\tt vladimir.gurvich@gmail.com}\\
\vskip 10pt

{\bf Mariya Naumova}\\
{\smallit Business School, Rutgers University,  Piscataway, NJ, United States}\\
{\tt mnaumova@business.rutgers.edu}\\

\end{center}

\vskip 20pt

\centerline{\bf Abstract}
Given integers $m > 1$ and $n > 0$,  
we say that a partition of the set $[n] = \{1, \dots, n\}$ is {\em $m$-good} 
if the number of elements in each part is at most  $m$ 
and their sum is a power of  $m$. 
It is easily seen that for every  $n$  there is a unique 2-good partition of $[n]$ and 
for each $m > 3$ there is no $m$-good partition for infinitely many $n$.   
Less is known for $m=3$. 
\newline
We conjecture that a 3-good partition of $[n]$ exists for each $n$ and 
prove that a minimal counter-example, if any, must be of the form: 
\newline 
(i) $n = 3^t + 3k +2$, where  $t > 0$  and 
(ii) $0 \leq k < \frac{3^{t-1}-1}{2}$; 
moreover, (iii) $k \neq \frac{3^\ell - 1}{2}$   
for all nonnegative integers $\ell < t$. 

Obviously, these conditions can be equivalently rewritten as: 
\newline 
(i$'$) $n \equiv 2 \; \pmod 3$, 
(ii$'$) $3k + 2 < \frac{3^{t+1} + 1}{2}$, and 
(iii$'$) $3k + 2 \neq \frac{3^{\ell + 1} + 1}{2}$   for $0 \leq \ell < t$. 
\newline 

By computations, the above conjecture was verified  for $n \leq 844$. 
\newline 
We also modify the statement slightly and 
prove it for the 3-good quasi-partitions, 
which cover all numbers 
of $[n] = \{1, \dots, 3^t+3k+2\}$ once, 
except  $3^t$, which is covered twice.
\newline 
Finally, we prove that a 3-good partition of $[n]$ is unique if 
\newline 
{\centering $n \in \{1,2,3,4, 3^t-4, 3^t-2, 3^t-1, 3^t, 
3^t+1, 3^t+2, 3^t+3, 3^t+5 \;\; \text{for} \;\; 
t \geq 2\}$}, 

\noindent 
and there are exactly two 3-good partitions of $[n]$ for $n = 3^t-3$. 
We conjecture that the number of 3-good partitions 
is greater than 2  for any other $n$, except 13.


\section{Introduction}
A partition of a positive integer $n$ is defined in \cite{Ald69} 
as a way of representing  $n$ as the sum of positive integers. 
Surveys on partitions can be found in \cite{Ald69}, \cite{AE04}.

For example, partitions of 3 can be expressed as (3), (1+2), and (1+1+1). 

Given integer $n > 0$ and $m > 1$, 
we will say that a subset of $[n] = \{1, \dots, n\}$  is {\em $m$-good} 
if its size is at most $m$ and 
the sum of numbers in it is a power of $m$, that is, $m^t$ for some integer $t \geq 0$. 
In particular, singleton $\{1\}$ is $m$-good for every $m$  and  $n$.

Furthermore, a partition of  $[n] = \{1, \dots, n\}$ will be called {\em $m$-good} 
if each of its partitioning sets (called {\em parts}) F is $m$-good. 

For example, $\{\{1,2,6\},\{3\},\{4,5\}\}$ is a 3-good partition of [6]. 
To  simplify notation, we will write $21 = (1+2+6) + (3) + (4+5)$ instead.
The number $T_n = 1 + \dots + n = n(n+1)/2$  is called {\em triangular}. 
We will denote a partition of $[n]$ 
by representing $T_n$ as the sum of the corresponding clauses.   
The summands within a sum are written in the increasing order, 
while the sums are ordered according to their smallest summands. 
As a brief reminder, the value on the left side of each equation below is the triangular number $T_n$, representing the total sum of the elements in the partitioned set $[n]$.

The first nine and seven $m$-good partitions for $m=2$ and $3$, respectively, are as follows:

$$1=(1), \; 3 = (1)+(2), \;6 = (1+3)+(2), \; 10 = (1+3)+(2)+(4), \; 
15 = (1)+(2)+(3+5)+(4),$$ 
$$21 = (1)+(2+6)+(3+5)+(4), \; 
28=(1+7)+(2+6)+(3+5)+(4), \;$$    
$$36 = (1+7)+(2+6)+(3+5)+(4)+(8), \; 
45=(1)+(2+6)+(3+5)+(4)+(8)+(7+9), \ldots$$    

$$1=(1), \; 3 = (1+2), \;6 = (1+2)+(3), \; 10 = (1)+(2+3+4), \; 
15 = (1+2)+(3+5)+(4),$$ 
$$21 = (1+2)+(3+6)+(4+5)=(1+2+6)+(3)+(4+5), \; 
28=(1)+(2+7)+(3+6)+(4+5), \ldots$$    

\smallskip 

Note that there exist two distinct $3$-good partitions for $n=6$. 

Note also that, obviously, there are no $m$-good partitions of $[n]$ for $n=2$  and $m > 3$. 

The following stronger claim shows that for any fixed $m>3$ 
there exists no $m$-good partition for infinitely many  $n$. 

\begin{observation}
For integer $m > 3$, $t \geq 0$, and $n = 2(mt+1)$, there is no $m$-good partition of $[n]$. 
\end{observation}

 \proof 
 Indeed,  if $m>3$ then $T_n = n(n+1)/2 = (mt+1)(2mt+3) \equiv 3\pmod m$, 
 while it must be $0$ or $1\pmod m$, by the definition of $m$-good partitions. 
 \qed 

\medskip 

It is also easily seen that the triangular numbers $T_n$ can only be congruent to $0$ or $1$ modulo $3$; that is, $T_n \not\equiv 2 \pmod{3}$.

\medskip 

\begin{proposition}
A unique $2$-good partition exists for any $n$. 
\end{proposition}

\proof 
Without loss of generality, 
assume that 
$n = 2^t + i$, where $t \geq 0$  and $0 \leq i < 2^t$. 
It is easy to verify that, for any $j$ such that $1 \leq j \leq i$, 
number $2^t+j$ can form a pair only with $2^t-j$, 
in any 2-good partition of $[n]$. 
In case $j=0$, both numbers equal $2^t$ and, hence,  
$2^t$ must participate in a 2-good partition of $[n]$ as a singleton. 
In any case, in $i+1$ steps, the problem of 2-good partitioning of $[n]$ 
is reduced to the same problem for $[n-2i-1]$, 
which is a strict subset of $[n]$. 
Applying this reduction recursively, we obtain 
a unique 2-good partition of $[n]$.

\qed 

\medskip

From now on we will consider only the case $m=3$, which looks difficult.  

\begin{conjecture}
\label{cnj-main}
A $3$-good partition of $[n]$ exists for every positive integer $n$.
\end{conjecture}

We will prove the following partial results: 

\begin{theorem}
\label{t-main}
A minimal counter-example, if any, must be of the form: 
\newline 
(i) $n = 3^t + 3k +2$, where  $t > 0$  and (ii) $0 \leq k < \frac{3^{t-1}-1}{2}$; 
\newline 
moreover, (iii) $k \neq \frac{(3^\ell - 1)}{2}$   
for all nonnegative integers $\ell < t$. 
\end{theorem} 

More partial results can be found in \cite{GN25}. 

We make several trivial preliminary observations regarding this statement:   

Obviously, (i) implies that  $n \equiv 2 \pmod 3$,   
and (ii, iii)  can be rewritten as: 
\newline 
(ii$'$) $3k + 2 < \frac{3^{t+1} + 1}{2}$ and 
(iii$'$) $3k + 2 \neq \frac{3^{\ell + 1} + 1}{2}$   for $0 \leq \ell < t$. 

Note that for $\ell > 0$, the values $\frac{3^\ell + 1}{2}, \; \frac{3^\ell - 1}{2}$, and 
$\frac{3^\ell - 1}{2} - 1 = 3 \cdot \frac{3^{\ell - 1} - 1}{2}$  
are three successive non-negative decreasing integers 
congruent to 2, 1, and 0 modulo 3, respectively.

\smallskip 

Statement (iii) excludes  $k = 0,1,4,13,40,\dots$   
and, respectively, $3k+2 = 2,5,14,41,122,\dots$  

Thus, the first numbers not covered by Theorem \ref{t-main} are 
$n = 35, 38, 89, 92, 98, 101, \dots$   

\medskip 

We will also prove the following related result.
A {\em quasi-partition} of $[n]$   
is a family of its subsets  
that covers one fixed number twice, while covering the remaining $n-1$ numbers exactly once each.  
Furthermore, we call a quasi-partition {\em $m$-good} 
if every one of its subsets is of cardinality at most  $m$  
and the sum of its numbers is a power of $m$, as before. 

\begin{theorem}
\label{t-quasi}
For each $n = 3^t + 3k + 2$  such that 
$t > 0$  and  $0 \leq k < \frac{3^{t-1}-1}{2}$  
there exists a $3$-good quasi-partition 
of $[n]$  in which $3^t$  is covered twice. 
\end{theorem}

\medskip 

We also obtain the following results on the uniqueness 
of 3-good partitions. 

\begin{theorem}
\label{t-unique}
A 3-good partition of $[n]$ is unique for  

{\centering $n \in N_u = \{1,2,3,4, 3^t-4, 3^t-2, 3^t-1, 3^t, 
3^t+1, 3^t+2, 3^t+3, 3^t+5 \;\; \text{for} \;\; 
t \geq 2\}$}. 

\noindent 
Furthermore, there are exactly two 3-good partitions of $[n]$ for $n = 3^t-3$ 
and also for $n=13$. 
\end{theorem}

We conjecture that the number of 3-good partitions 
is greater than 2  for any other $n$. 

This conjecture and Conjecture \ref{t-main},  
were verified up to $n = 844$  
by a computer code written by Dmitry Rybin \cite{Ryb25}. 
The next observation is due to Oleg Polubasoff. 

\begin{proposition} (\cite{Pol25}) 
\label{p-Polubas}
A 3-good partition of  $[n]$  is uniquely defined 
by the set of its triplets (and by $n$).  
\end{proposition}

\section{Proof of Theorem \ref{t-main}} 
Denote by $\operatorname{part}(n)$  a 3-good partition of $[n]$.
(The reader should not confuse it with partition function $p(n)$.) 
Assume for contradiction that there exists an $n$  without $\operatorname{part}(n)$. 
Choose the minimal such $n$ and call it {\em critical}. 

\subsection{Reduction}
\label{ss11}
Obviously, $n$ is not critical if  
$\operatorname{part}(n) = \operatorname{part}(n') + (C_1) + \ldots + (C_\ell)$, 
where all these $\ell$ added parts are of size at most 3 and partition the interval $[n'+1,n]$.
Such a representation of  $\operatorname{part}(n)$  will be called a reduction of  $n$ to $n'$ 
or just a reduction of $n$, for short.  

For example, $n = 3^t$ and $3^t+1$ admit reductions  
$\operatorname{part}(3^t) = \operatorname{part}(3^t-1) +(3^t)$  and 
$\operatorname{part}(3^t+1) = \operatorname{part}(3^t-2) + ((3^t-1) + 3^t + (3^t+1))$. Also, 

\begin{align*}
\operatorname{part}(3) &= \operatorname{part}(2) + (3), \\
\operatorname{part}(4) &= \operatorname{part}(1) + (2 + 3 + 4), \\
\operatorname{part}(5) &= \operatorname{part}(1) + (2 + 3 + 4), \\
\operatorname{part}(6) &= \operatorname{part}(2) + (3 + 6) + (4 + 5), \\
\operatorname{part}(7) &= \operatorname{part}(1) + (2 + 7) + (3 + 6) + (4 + 5), \\
\operatorname{part}(9) &= \operatorname{part}(8) + (9), \\
\operatorname{part}(10) &= \operatorname{part}(7) + (8 + 9 + 10).
\end{align*}

\smallskip

Note that when \( n = 3^t - 1 \), the value is not reducible in the usual sense, but it admits a unique 3-good partition:
\[
\operatorname{part}(3^t - 1) = (1 + (3^t - 1)) + (2 + (3^t - 2)) + \cdots + \left( \frac{3^t - 1}{2} + \frac{3^t + 1}{2} \right).
\]
\noindent 
Note that here, as well as in some further similar formulas, 
clauses may denote numbers as well as singleton parts. 
In the above formula, each clause contains two numbers.

For example:
\[
\operatorname{part}(2) = (1 + 2), \quad 
\operatorname{part}(8) = (1 + 8) + (2 + 7) + (3 + 6) + (4 + 5), \quad \dots
\]

\smallskip

The next non-reducible value is \( n = 11 \), which has a more intricate 3-good partition:
\[
\operatorname{part}(11) = (1 + 8) + (2 + 7) + (3) + (4 + 5) + (6 + 10 + 11) + (9).
\]

\smallskip

In general, for \( n = 3^t + 2 \), we have:
\[
\operatorname{part}(3^t + 2) = (1 + (3^t - 1)) + (2 + (3^t - 2)) + (3) + (4 + (3^t - 4)) + 
\cdots + \left( \frac{3^t - 1}{2} + \frac{3^t + 1}{2} \right)
\]
\[
\hspace{9mm} +\; ((3^t - 3) + (3^t + 1) + (3^t + 2)) + (3^t).
\]

\subsection{Proof of Theorem \ref{t-main} (ii)}
\label{ss12}
Clearly, for each  $n > 0$  there is a unique $t > 0$ such that 
$3^{t-1} \leq n < 3^t$.   
Theorem \ref{t-main}~(ii) asserts that inequalities $3^{t-1} \leq n \leq \frac{3^t-1}{2}$ must hold whenever $n$ is critical.  
  
For example, consider the following reductions:
\begin{align*}
\operatorname{part}(14) &= \operatorname{part}(12) + (13 + 14), \\
\operatorname{part}(15) &= \operatorname{part}(11) + (12 + 15) + (13 + 14), \\
\operatorname{part}(25) &= \operatorname{part}(1) + (2 + 25) + (3 + 24) + (4 + 23) + \dots + (12 + 15) + (13 + 14).
\end{align*}

\smallskip

\noindent
Note that $15 \equiv 0 \pmod{3}$ is reduced to $11 \equiv 2 \pmod{3}$.

\smallskip

\noindent
In general, in accordance with Theorem \ref{t-main} (ii), any integer \( n \) satisfying
\[
\frac{3^t + 1}{2} \leq n < 3^t,
\]
is reducible. Indeed, the reduction is given by the following formula:
\[
\operatorname{part}\left(\frac{3^t + 1}{2} + \ell\right) = 
\operatorname{part}\left(\frac{3^t - 1}{2} - \ell - 1\right) + 
\sum_{i=0}^{\ell} \left(\left(\frac{3^t - 1}{2} - i\right) + \left(\frac{3^t + 1}{2} + i\right)\right),
\]
where \( t \geq 0 \) and \( 0 \leq \ell < \frac{ 3^t + 1}{2} \). This establishes the claim.

\medskip

\noindent
Note that no new triples appear in the reduction; in other words, the 3-good partitions
\[
\operatorname{part}\left(\frac{3^t + 1}{2} + \ell\right) \quad \text{and} \quad \operatorname{part}\left(\frac{3^t - 1}{2} - \ell - 1\right)
\]
contain the same set of triples.

\smallskip

\noindent
The triple-free 3-good partitions \( \operatorname{part}(3^t) \) and \( \operatorname{part}(3^t - 1) \) 
were discussed in Section~\ref{ss11}. 

\subsection{3-good partitions of symmetric integer intervals}
\label{ss13}
Let  $[-\ell;\ell] = \{-\ell,\dots,-1,0,1,\dots,\ell\}$ be the symmetric integer interval of size $2\ell+1$. 

For  $\ell \equiv 1 \pmod 3$, we call a partition of $[-\ell;\ell]$, 
denoted by $\operatorname{part}[-\ell;\ell]$, {\em 3-good} if every of its parts is of size 3 and the sum of the numbers in it equals to 0. 
For  $\ell \equiv 0 \pmod 3$, we slightly modify this definition by adding one special part 
$\{0\}$ of size $1$, containing only $0$.   

\smallskip 

Note that for  $\ell \equiv 2 \pmod 3$  the size of $[-\ell;\ell]$ is $2\ell+1 \equiv 2 \pmod 3$. 
In this case, we have no 3-good partitions. 

\begin{lemma}
If  $\ell \equiv 0$ or $1 \pmod 3$ then $[-\ell;\ell]$ admits a 3-good partition.  
\end{lemma}

\proof  
Assume first that \( \ell = 3k \). Then, the partition \( \operatorname{part}[-3k; 3k] \) is given by:
\begin{align*}
\operatorname{part}[-3k; 3k] 
&= (0) + \sum_{i=1}^{k-1} \left( (-3k + i) + (1 + i) + (3k - 1 - 2i) \right) \\
&\quad + \sum_{i=1}^{k-1} \left( (-2k + i) + (-k + i) + (3k - 2i) \right).
\end{align*}

\noindent For example, if \( k = 1 \) and \( k = 2 \), we obtain:
\begin{align*}
\operatorname{part}[-3; 3] &= (-3 + 1 + 2) + (-2 + (-1) + 3) + (0), \\
\operatorname{part}[-6; 6] &= (-6 + 1 + 5) + (-5 + 2 + 3) + (-4 + (-2) + 6) \\
           &\quad + (-3 + (-1) + 4) + (0).
\end{align*}

\medskip 

For  $\ell = 3k+1$  replace  $(0)$ by $(-\ell+0+\ell)$.
\qed

\medskip 

However, $n \equiv 2 \pmod 3$ are not reducible in this way. 
Indeed, in this case, interval $[-n;n]$ admits no 3-good partitions, since its size is $2n+1 \equiv 2 \pmod 3$. 
In particular, this holds for $n=3^t+3k+2$ for any $t > 0$ and $k \geq 0$. 

\subsection{Proof of Theorem \ref{t-main} (i)} 
\label{ss14}
Let us show that $n \equiv 2 \pmod 3$ must hold for each critical $n$. 
If  $n \equiv 0$ or $1 \pmod 3$ then reductions work.
Consider, for example, the following reductions:
\begin{align*}
\operatorname{part}(6) &= \operatorname{part}(2) + (3 + 6) + (4 + 5), \\
\operatorname{part}(12) &= \operatorname{part}(5) + (6 + 10 + 11) + (7 + 8 + 12) + (9), \\
\operatorname{part}(30) &= \operatorname{part}(23) + (24 + 28 + 29) + (25 + 26 + 30) + (27), \\
\operatorname{part}(33) &= \operatorname{part}(20) + (21 + 28 + 32) + (22 + 29 + 30) + (23 + 25 + 33) \\
&\quad + (24 + 26 + 31) + (27).
\end{align*}

More generally, for \( n = 3^t + 3 \), we have:
\[
\operatorname{part}(3^t + 3) = \operatorname{part}(3^t - 4) 
+ \left( (3^t - 3) + (3^t + 1) + (3^t + 2) \right) 
+ \left( (3^t - 2) + (3^t - 1) + (3^t + 3) \right).
\]

Eliminating terms involving \( 3^t \) yields the 3-good partition:
\[
\operatorname{part}[-3;3] = (-3 + 1 + 2) + (-2 - 1 + 3) + (0).
\]

\medskip 

\noindent
Consider also the following reductions:
\begin{align*}
\operatorname{part}(7)  &= \operatorname{part}(1) + (2 + 7) + (3 + 6) + (4 + 5), \\
\operatorname{part}(13) &= \operatorname{part}(4) + (5 + 9 + 13) + (6 + 10 + 11) + (7 + 8 + 12) + (9), \\
\operatorname{part}(31) &= \operatorname{part}(22) + (23 + 27 + 31) + (24 + 28 + 29) + (25 + 26 + 30), \\
\operatorname{part}(34) &= \operatorname{part}(19) + (20 + 30 + 31) + (21 + 28 + 32) + (22 + 26 + 33) + (23 + 24 + 34) + (25 + 27 + 29).
\end{align*}

\smallskip

\noindent
More generally, for \( n = 3^t + 4 \), we have:
\[
\operatorname{part}(3^t + 4) = \operatorname{part}(3^t - 5)
+ \big((3^t - 4) + 3^t + (3^t + 4)\big)
+ \big((3^t - 3) + (3^t + 1) + (3^t + 2)\big)
+ \big((3^t - 2) + (3^t - 1) + (3^t + 3)\big).
\]

\noindent
Eliminating all instances of \( 3^t \) yields the 3-good partition:
\[
\operatorname{part}[-4;4] = (-4 + 0 + 4) + (-3 + 1 + 2) + (-2 -1 + 3).
\]

\medskip

\begin{figure}[htbp]
\begin{center}
\begin{tikzpicture}[scale=1.5]
\def\colwidth{1}
\def\rowheight{0.6}
\def\nrows{3}
\def\gap{0.4}

\def\xLeft{0}
\def\xRight{2*\colwidth + 2*\gap}
\def\xMid{\colwidth + \gap}

\foreach \i in {0,1,2,3} {
    \draw (\xLeft,-\i*\rowheight) -- (\xLeft+\colwidth,-\i*\rowheight);
}
\draw (\xLeft,0) -- (\xLeft,-\nrows*\rowheight);
\draw (\xLeft+\colwidth,0) -- (\xLeft+\colwidth,-\nrows*\rowheight);

\foreach \i in {0,1,2,3} {
    \draw (\xRight,-\i*\rowheight) -- (\xRight+\colwidth,-\i*\rowheight);
}
\draw (\xRight,0) -- (\xRight,-\nrows*\rowheight);
\draw (\xRight+\colwidth,0) -- (\xRight+\colwidth,-\nrows*\rowheight);

\def\yRect{0}
\draw (\xMid,\yRect) rectangle ++(\colwidth,\rowheight);

\def\labelsize{\scriptsize}

\node at (\xLeft + 0.5*\colwidth,-0.5*\rowheight) {\small \( 3^t - 1 \)};
\node at (\xLeft + 0.5*\colwidth,-1.5*\rowheight) {\small \( 3^t - 2 \)};
\node at (\xLeft + 0.5*\colwidth,-2.5*\rowheight) {\small \( 3^t - 3 \)};

\node at (\xRight + 0.5*\colwidth,-0.5*\rowheight) {\small \( 3^t + 1 \)};
\node at (\xRight + 0.5*\colwidth,-1.5*\rowheight) {\small \( 3^t + 2 \)};
\node at (\xRight + 0.5*\colwidth,-2.5*\rowheight) {\small \( 3^t + 3 \)};

\node[anchor=center] at 
  (\xMid + 0.5*\colwidth, \yRect - 0.5*\rowheight + 0.65) 
  {\small \( 3^t \)};


\draw (\xLeft + \colwidth, -0.5*\rowheight) -- (\xRight, -2.5*\rowheight);

\draw (\xLeft + \colwidth, -1.5*\rowheight) -- (\xRight, -2.5*\rowheight);

\draw (\xLeft + \colwidth, -2.5*\rowheight) -- (\xRight, -0.5*\rowheight);

\draw (\xLeft + \colwidth, -2.5*\rowheight) -- (\xRight, -1.5*\rowheight);

\end{tikzpicture}
\end{center}
\caption{Reduction for $n = 3^t + 3k$, $k=1, t \ge 2$} \label{f1}
\end{figure}
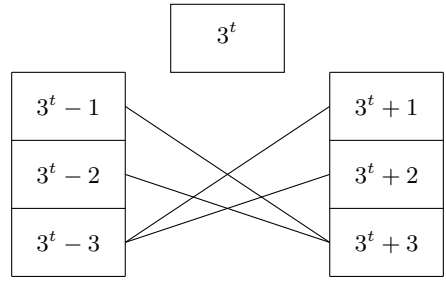

\begin{figure}[htbp]
\begin{center}
\begin{tikzpicture}[scale=1.5]

\def\colwidth{1}
\def\rowheight{0.6} 
\def\nrows{6}
\def\gap{0.4}

\def\xLeft{0}
\def\xRight{2*\colwidth + 2*\gap}
\def\xMid{\colwidth + \gap}

\foreach \i in {0,...,\nrows} {
    \draw (\xLeft,-\i*\rowheight) -- (\xLeft+\colwidth,-\i*\rowheight);
}
\draw (\xLeft,0) -- (\xLeft,-\nrows*\rowheight);
\draw (\xLeft+\colwidth,0) -- (\xLeft+\colwidth,-\nrows*\rowheight);

\foreach \i in {0,...,\nrows} {
    \draw (\xRight,-\i*\rowheight) -- (\xRight+\colwidth,-\i*\rowheight);
}
\draw (\xRight,0) -- (\xRight,-\nrows*\rowheight);
\draw (\xRight+\colwidth,0) -- (\xRight+\colwidth,-\nrows*\rowheight);

\def\yRect{0}
\draw (\xMid,\yRect) rectangle ++(\colwidth,\rowheight);

\def\labelsize{\scriptsize}

\foreach \i in {1,...,6} {
    \node at (\xLeft + 0.5*\colwidth, {-(\i-0.5)*\rowheight}) {\small \( 3^t - \i \)};
}

\foreach \i in {1,...,6} {
    \node at (\xRight + 0.5*\colwidth, {-(\i-0.5)*\rowheight}) {\small \( 3^t + \i \)};
}

\node[anchor=center] at 
  (\xMid + 0.5*\colwidth, \yRect - 0.5*\rowheight + 0.55) 
  {\small \( 3^t \)};



\draw (\xLeft + \colwidth, -0.5*\rowheight) -- (\xRight, -3.5*\rowheight);

\draw (\xLeft + \colwidth, -1.5*\rowheight) -- (\xRight, -5.5*\rowheight);

\draw (\xLeft + \colwidth, -2.5*\rowheight) -- (\xRight, -3.5*\rowheight);

\draw (\xLeft + \colwidth, -3.5*\rowheight) -- (\xRight, -5.5*\rowheight);

\draw (\xLeft + \colwidth, -4.5*\rowheight) -- (\xRight, -1.5*\rowheight);
\draw (\xLeft + \colwidth, -4.5*\rowheight) -- (\xRight, -2.5*\rowheight);

\draw (\xLeft + \colwidth, -5.5*\rowheight) -- (\xRight, -0.5*\rowheight);
\draw (\xLeft + \colwidth, -5.5*\rowheight) -- (\xRight, -4.5*\rowheight);

\end{tikzpicture}
\end{center}
\caption{Reduction for $n = 3^t + 3k$, $k=2, t \ge 3$} \label{f2}
\end{figure}
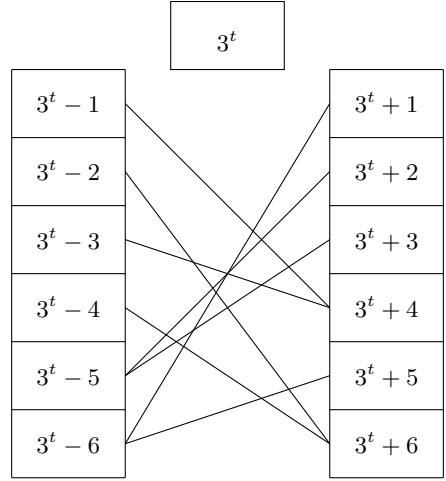

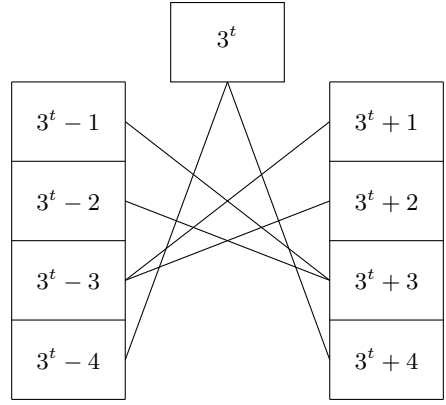
\begin{figure}[htbp]
\begin{center}
\begin{tikzpicture}[scale=1.5]
\def\colwidth{1}
\def\rowheight{0.7}
\def\nrows{4}
\def\gap{0.4}

\def\xLeft{0}
\def\xRight{2*\colwidth + 2*\gap}
\def\xMid{\colwidth + \gap}

\foreach \i in {0,...,\nrows} {
    \draw (\xLeft,-\i*\rowheight) -- (\xLeft+\colwidth,-\i*\rowheight);
}
\draw (\xLeft,0) -- (\xLeft,-\nrows*\rowheight);
\draw (\xLeft+\colwidth,0) -- (\xLeft+\colwidth,-\nrows*\rowheight);

\foreach \i in {0,...,\nrows} {
    \draw (\xRight,-\i*\rowheight) -- (\xRight+\colwidth,-\i*\rowheight);
}
\draw (\xRight,0) -- (\xRight,-\nrows*\rowheight);
\draw (\xRight+\colwidth,0) -- (\xRight+\colwidth,-\nrows*\rowheight);

\def\yRect{0}
\draw (\xMid,\yRect) rectangle ++(\colwidth,\rowheight);

\def\labelsize{\scriptsize}

\foreach \i in {1,...,4} {
    \node at (\xLeft + 0.5*\colwidth, {-(\i-0.5)*\rowheight}) {\small \( 3^t - \i \)};
}

\foreach \i in {1,...,4} {
    \node at (\xRight + 0.5*\colwidth, {-(\i-0.5)*\rowheight}) {\small \( 3^t + \i \)};
}

\node[anchor=center] at 
  (\xMid + 0.5*\colwidth, \yRect + 0.5*\rowheight + 0.05) 
  {\small \( 3^t \)};

\def\midBottomX{\xMid + 0.5*\colwidth}
\def\midBottomY{\yRect}


\draw (\xLeft + \colwidth, -0.5*\rowheight) -- (\xRight, -2.5*\rowheight);
\draw (\xLeft + \colwidth, -1.5*\rowheight) -- (\xRight, -2.5*\rowheight);
\draw (\xLeft + \colwidth, -2.5*\rowheight) -- (\xRight, -0.5*\rowheight);
\draw (\xLeft + \colwidth, -2.5*\rowheight) -- (\xRight, -1.5*\rowheight);

\draw (\xLeft + \colwidth, -3.5*\rowheight) -- (\midBottomX, \midBottomY);
\draw (\xRight, -3.5*\rowheight) -- (\midBottomX, \midBottomY);

\end{tikzpicture}
\end{center}
\caption{Reduction for $n = 3^t + 3k+1$, $k=1, t \ge 2$} \label{f3}
\end{figure}

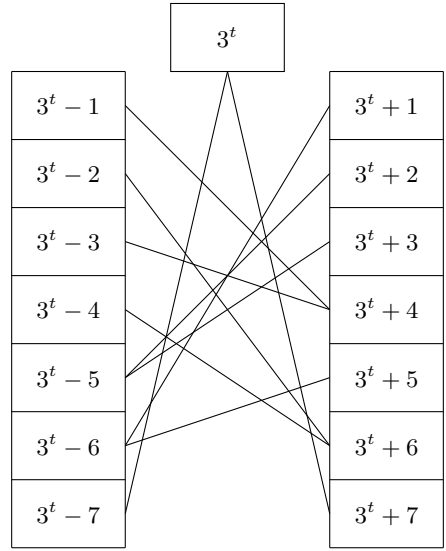
\begin{figure}[htbp]
\begin{center}
\begin{tikzpicture}[scale=1.5]
\def\colwidth{1}
\def\rowheight{0.6} 
\def\nrows{7}
\def\gap{0.4}

\def\xLeft{0}
\def\xRight{2*\colwidth + 2*\gap}
\def\xMid{\colwidth + \gap}

\foreach \i in {0,...,\nrows} {
    \draw (\xLeft,-\i*\rowheight) -- (\xLeft+\colwidth,-\i*\rowheight);
}
\draw (\xLeft,0) -- (\xLeft,-\nrows*\rowheight);
\draw (\xLeft+\colwidth,0) -- (\xLeft+\colwidth,-\nrows*\rowheight);

\foreach \i in {0,...,\nrows} {
    \draw (\xRight,-\i*\rowheight) -- (\xRight+\colwidth,-\i*\rowheight);
}
\draw (\xRight,0) -- (\xRight,-\nrows*\rowheight);
\draw (\xRight+\colwidth,0) -- (\xRight+\colwidth,-\nrows*\rowheight);

\def\yRect{0}
\draw (\xMid,\yRect) rectangle ++(\colwidth,\rowheight);

\def\labelsize{\scriptsize}

\foreach \i in {1,...,7} {
    \node at (\xLeft + 0.5*\colwidth, {-(\i-0.5)*\rowheight}) {\small \( 3^t - \i \)};
}

\foreach \i in {1,...,7} {
    \node at (\xRight + 0.5*\colwidth, {-(\i-0.5)*\rowheight}) {\small \( 3^t + \i \)};
}

\node[anchor=center] at 
  (\xMid + 0.5*\colwidth, \yRect - 0.5*\rowheight + 0.6) 
  {\small \( 3^t \)};



\draw (\xLeft + \colwidth, -0.5*\rowheight) -- (\xRight, -3.5*\rowheight); 
\draw (\xLeft + \colwidth, -1.5*\rowheight) -- (\xRight, -5.5*\rowheight); 
\draw (\xLeft + \colwidth, -2.5*\rowheight) -- (\xRight, -3.5*\rowheight); 
\draw (\xLeft + \colwidth, -3.5*\rowheight) -- (\xRight, -5.5*\rowheight); 
\draw (\xLeft + \colwidth, -4.5*\rowheight) -- (\xRight, -1.5*\rowheight); 
\draw (\xLeft + \colwidth, -4.5*\rowheight) -- (\xRight, -2.5*\rowheight); 
\draw (\xLeft + \colwidth, -5.5*\rowheight) -- (\xRight, -0.5*\rowheight); 
\draw (\xLeft + \colwidth, -5.5*\rowheight) -- (\xRight, -4.5*\rowheight); 

\draw (\xLeft + \colwidth, -6.5*\rowheight) -- (\xMid + 0.5*\colwidth, \yRect); 
\draw (\xRight, -6.5*\rowheight) -- (\xMid + 0.5*\colwidth, \yRect);          

\end{tikzpicture}
\end{center}
\caption{Reduction for $n = 3^t + 3k+1$, $k=2, t \ge 3$} \label{f4}
\end{figure}

\noindent
In general, the reductions for \( n = 3^t + 3k \) and \( n = 3^t + 3k + 1 \) 
are given by the following formulas (see also examples in Figures \ref{f1}, \ref{f2}, \ref{f3}, \ref{f4}):

\smallskip

\begin{align*}
\operatorname{part}(3^t + 3k) &= \operatorname{part}(3^t - 3k - 1) + (3^t) \\
&\quad + \sum_{i = 0}^{k - 1} \left((3^t - 3k + i) + (3^t + i + 1) + (3^t + 3k - 2i - 1)\right) \\
&\quad + \sum_{i = 0}^{k - 1} \left((3^t - 2k + i) + (3^t - k + i) + (3^t + 3k - 2i)\right).
\end{align*}

\smallskip

\begin{align*}
\operatorname{part}(3^t + 3k + 1) &= \operatorname{part}(3^t - 3k - 1) \\
&\quad + \left((3^t - 3k - 1) + 3^t + (3^t + 3k + 1)\right) \\
&\quad + \sum_{i = 0}^{k - 1} \left((3^t - 3k + i) + (3^t + i + 1) + (3^t + 3k - 2i - 1)\right) \\
&\quad + \sum_{i = 0}^{k - 1} \left((3^t - 2k + i) + (3^t - k + i) + (3^t + 3k - 2i)\right).
\end{align*}

\smallskip

\vspace{2cm} 

\noindent
In both cases, the condition \( 3^{t-1} \geq 2k + 1 \) must hold, which is equivalent to
\[
k \leq \frac{3^{t-1} - 1}{2}.
\]

\noindent
Recall that if \( 3^t + 1 \geq k > \frac{3^{t-1} - 1}{2} \), 
then the reduction described in Section~\ref{ss12} applies.

\smallskip 

Thus, the two reductions discussed above apply to all integers \( n \), 
except those satisfying both of the following conditions:
\[
3^t + 1 < n < \frac{3^{t+1} + 1}{2} \quad \text{and} \quad n = 3^t + 3k + 2,
\]
for some nonnegative integers \( t \) and \( k \).

\smallskip

Note that the case of \( n \equiv 1 \pmod{3} \) is covered by the reduction, 
while case of \( n \equiv 0 \pmod{3} \) is reduced to \( n' \equiv 2 \pmod{3} \).

\subsection{Proof of Theorem \ref{t-main} (iii)} 
\label{ss(iii)}
Recall that for $\ell > 0$, the values $\frac{3^\ell + 1}{2}, \; \frac{3^\ell - 1}{2}$, and 
$\frac{3^\ell - 1}{2} - 1 = 3 \cdot \frac{3^{\ell - 1} - 1}{2}$  
are three successive non-negative decreasing integers 
congruent to 2, 1, and 0 modulo 3, respectively.
Hence, $n = 3^t + \frac{3^\ell + 1}{2} \equiv 2 (\pmod{3})$. 
Consider triplet 
$T = ((3^t + \frac{3^\ell + 1}{2}) + (3^t + \frac{3^\ell - 1}{2}) + (3^t - 3^\ell))$ 
whose sum, $3^{t+1}$, is a power of  3.  

Suppose for contradiction that $n$ is critical. 
Then there is a 3-good partition of $[n-2]$. 
Obviously, one of its  sets $T'$ contains the number $(3^t - 3^\ell)$. 
Then, $T'$ contains one or two more numbers. Clearly, their sum is $3^\ell$. 
Hence, we can just move $(3^t - 3^\ell)$  from  $T'$ to $T$ 
and obtain a 3-good partition of  $[n]$,   
together with a reduction: $n$ to $n-2$. 
Thus, $n$ is not critical, which is a contradiction. 

\smallskip 

In \cite{GN25} several subcases left open by Theorem \ref{t-main} 
were resolved 
by a more powerful reduction suggested by Oleg Polubasoff. 

\section{Proof of Theorem \ref{t-quasi}} 
Fix a nonnegative integer $k$ and set $n = 3k+2$.  
We will partition $[n]$  into  $k$  triples 
$T_i = (a_i, b_i, c_i), i = 1, \dots, k$ 
such that  $a_i = b_i + c_i$  and two singletons: 
$U = 3k+2$  and  $L = \lceil{(3k+1)/2}\rceil$.

Define the triples as follows: 
$c_i$  is a sequence of length  $k$   
in which the numbers that are not multiples of 3 
go in the monotone increasing order and 
then the numbers multiple of 3  
go in the monotone decreasing order. 
For example, for  $k = 0,1,\dots,11$   
we obtain: 

\smallskip 

\noindent 
$\emptyset; \; 1; \; 1,2; \; 1,2,3; \; 1,2,4,3; \; 1,2,4,5,3; 
\; 1,2,4,5,6,3; \; 1,2,4,5,7,6,3; \; 1,2,4,5,7,8,6,3; 
\newline 
\; 1,2,4,5,7,8,9,6,3; \; 1,2,4,5,7,8,10,9,6,3; \; \; 1,2,4,5,7,8,10,11,9,6,3; \;  \dots$ 

Thus, we have  $c_i$; now to define  $a_i$  and  $b_i$ for $i = 1, \dots k$  
we proceed recursively as follows.
Delete $U, L$, and $[k]$  from $[n]$. 
For  $i = 1, \dots k$  define  $a_i$  as the maximum of the remaining numbers 
and set  $b_i = a_i - c_i$. Then proceed with $i+1$.  
For example, for  $k = 0,1,\dots,7$  we obtain the following 
$U = U_k, L=L_k$, and $(a_i,b_i,c_i)$ for $i \in [k]$: 

$k=0: \;\; 2, 1$; 

$k=1: \;\; 5, 2; \; (4,3,1)$; 

$k=2: \;\; 8, 4; \; (7,6,1), (5,3,2)$;   

$k=3: \;\; 11, 5; \; (10,9,1), (8,6,2), (7,4,3)$;   

$k=4: \;\; 14, 7; \; (13,12,1), (11,9,2), (10,6,4), (8,5,3)$;   

$k=5: \;\; 17, 8; \; (16,15,1), (14,12,2), (13,9,4), (11,6,5), (10,7,3)$; 

$k=6: \;\; 20, 10; \; (19,18,1), (17,15,2), (16,12,4), (14,9,5), (13,7,6), (11,8,3)$; 

$k=7: \;\; 23, 11; \; (22,21,1), (20,18,2), (19,15,4), (17,12,5), (16,9,7), (14,8,6), (13,10,3)$.  

\medskip 

It is not difficult to verify that the desired partitions are well defined for all $k$. 

\medskip 

Replace each triple $(a_i,b_i,c_i)$ by two zero-sum triples 
$(a_i,-b_i,-c_i)$  and  $(-a_i,b_i,c_i)$; 
then, replace $U$ and $L$ by the zero-sum triples   
$\{U,0,-U\}$ and  $\{L,0,-L\}$, respectively. 

We obtain $2(k+1)$ zero-sum triples that cover 
each of $6k+5$ numbers of the symmetric integer interval $[-(3k+2); 3k+2]$ 
exactly once, except 0, which is covered twice. 

The above two constructions are of independent interest. 

Furthermore, increase each number by $3^t$ and   
add $part((3^t-3k-3))$ to the obtained set of triples. 
Note that such a partition exists, 
since $(3^t-3k-3) = 3 (3^{t-1}-k-1) \equiv 0 \; \pmod 3$.
Obviously, the above construction results in a 3-good quasi-partition 
of $[3^t + 3k + 2]$  by  $ 3^{t-1}+k+1$ triples in which number $3^t$  appears twice, 
in triples $(3^t - L, 3^t, 3^t + L)$  and $(3^t - U, 3^t, 3^t + U)$ .  
This proves Theorem~\ref{t-quasi}. 

\section{Proof of Theorem \ref{t-unique} and Proposition \ref{p-Polubas}}  
{\bf Case} $n = 3^t-4$. A unique 3-good partition is as follows: 

$$(1) + (2 + 3) + (4 + (3^t-4)) + (5 + (3^t-5)) \dots + 
\left( \frac{3^t - 1}{2} + \frac{3^t + 1}{2} \right).$$
 
Indeed, $3^t-4$  can be either in the pair with 4, or in the triple  with 1 and 3. 
Yet, the second option is impossible because of $3^t-3$.  

\medskip 

{\bf Case} $n = 3^t-3$.  There exist exactly two 3-good partitions:  

$(1 + 2) + (3 + (3^t-3)) + (4 + (3^t-4)) +  (5 + (3^t-5)) + 
\dots + \left( \frac{3^t - 1}{2} + \frac{3^t + 1}{2} \right)$
and 

$(1 + 2 +  (3^t-3)) + (3) +    
(4 + (3^t-4)) + (5 + (3^t-5)) + \dots +  \left( \frac{3^t - 1}{2} + \frac{3^t + 1}{2} \right)$.  

Indeed, there are only the above two ways to include $3^t-3$

\medskip 

{\bf Case} $n = 3^t-2$. A unique 3-good partition is as follows: 

$(1) + (2 + (3^t-2)) + (3 + (3^t-3)) +  \dots +  
\left( \frac{3^t - 1}{2} + \frac{3^t + 1}{2} \right)$.  

This can be proven by induction: 
$3^t-2$ can be only paired with 2; 
then, $3^t-3$ can be only paired with 3, 
since 2 is already taken, etc.; finally, 
1  remains as a legal singleton. 

\medskip 

{\bf Case}  $n = 3^t-1$. A unique 3-good partition consists of pairs only: 

$(1 + (3^t-1)) + (2 + (3^t-2)) +  \dots + 
\left( \frac{3^t - 1}{2} + \frac{3^t + 1}{2} \right)$. 

Again induction: 
$3^t-1$ can be only paired with 1  and then, in general, 
$3^t-k$  can be only paired with $k$    
for $1 \leq k \leq (3^t-1)/2$. 

\medskip 

{\bf Case} $n = 3^t$. A unique 3-good partition exists. 
Obviously, $(3^t)$ can appear only as a singleton and,  
by this, the considered case is reduced to the previous one. 

\medskip 

{\bf Case} $n = 3^t+1$. A unique 3-good partition exists. 
Obviously, triple
$(3^t-1, 3^t, 3^t+1)$  must appear, 
since this is the only way to include $3^t+1$. 
By this, the considered case is reduced to $n'= 3^t-2$.  

\medskip 

{\bf Case}  $n = 3^t+2$. A unique 3-good partition exists: 

$((3^t+2) + (3^t+1) + (3^t-3)) + 
(3^t) +$

$(1 + (3^t-1)) + (2 + (3^t-2)) + (3) +  
(4 + (3^t-4)) + \dots +  
\left( \frac{3^t - 1}{2} + \frac{3^t + 1}{2} \right)$. 

It is not difficult to verify that 
this is the only way to employ $3^t+2$. 
Indeed, the only alternative option 
$(3^t+2, 3^t, 3^t-2)$  is a dead-end, 
since after this including $3^t+1$ requires $3^t-1$ 
and  $3^t$, which is already taken. 

Note that the issue will be resolved if 
we are allowed to include $3^t$ twice. 

\medskip 

{\bf Case} $n = 3^t+3$. A unique 3-good partition exists: 

$((3^t+3) + (3^t-1) + (3^t-2)) + 
((3^t+2) + (3^t+1) + (3^t-3)) +$ 

$(3^t) + (4 + (3^t-4)) + \dots + \left( \frac{3^t - 1}{2} + \frac{3^t + 1}{2} \right).$ 

It is not difficult to verify that 
this is the only way to include $3^t+3$. 
Indeed, the only alternative option 
$((3^t+3) + 3^t + (3^t-3))$  is a dead-end, 
since after this, the inclusion $3^t+1$ 
(respectively, $3^t+2$)  requires $3^t-1$ 
(respectively, $3^t-2$)
and  $3^t$, which is already taken.
Again, the issue will be resolved if we 
are allowed to include $3^t$ twice. 

\medskip 

{\bf Case}  $n = 3^t+4$. 
It is easy to verify that $n=13$ 
admits exactly two 3-good partitions: 

$(13+8+6)+(12+10+5)+(11+9+7)+(4+3+2)+(1)$  and 

$(13+9+5)+(12+8+7)+(11+10+6)+(4+3+2)+(1)$.  

Furthermore,  $n=31$  admits fourteen 3-good partitions, 
according to our computations. 

\medskip  

{\bf Case} $n = 3^t+5$. A unique 3-good partition exists. 
It contains the triple 
\newline
$((3^t+5) +  (3^t+4) + (3^t-9))$.  
Furthermore, a unique 3-good partition 
for $n'=3^t+3$  contains the pair $(9 + (3^t-9))$. 
Although  $3^t-9$ is already taken, the remaining number $9=3^2$ is a proper singleton.
It is not difficult to verify that 
there is no other way to include $3^t+5$ and  $3^t+4$, since   
for every other option number $3^t$  appears more than once. 

\medskip 

Note that 3-good partitions $part(3^t - 1)$  and  $part(3^t)$ 
are unique, do not contain triples, and 
contain 0 and 1 singleton, respectively.

\bigskip 

Proposition \ref{p-Polubas} 
results from the following recursive algorithm 
suggested by Oleg Polubasoff \cite{Pol25}. 
A number from $[n]$ is called {\em free}  
if it is not yet included into a set of the partition under construction.   
Consider the largest free number  $i$. 
If  $i$  is a power of 3 then it is a singleton.  
Otherwise,  $3^{j-1} < i < 3^j$  for some $j$. 
In this case $i$  must be in pair with $3^j – i$.
If this number is already taken 
(that is, contained in a triple) or 
is larger than  $n$  then 
there exists no 3-good partition of $[n]$ 
with the given set of triples. 
If all numbers are taken, we obtain a unique desired partition.
For example, let the input be $n=4$  with no triples. 
The largest free number is 4. We have $3 < 4 < 9$, but $9-4=5>4=n$.  
Hence, there is no pair for 4 and, indeed, 
[4]  admits a unique 3-good partition (1) + (2+3+4), 
which contains a triple. 

\section{Generalized reductions}
Oleg Polubasoff \cite{Pol25} suggested a more general method.  
To reduce  $\operatorname{part}(n)$  to  $\operatorname{part}(n')$, where $n > n'$, 
we are allowed to delete some sets from $\operatorname{part}(n')$ and then 
add some sets to obtain $\operatorname{part}(n)$. 
This method can be successfully applied to the smallest open cases 
$n = 3^t + 8$  and  $3^t + 11$  for $t \geq 3$ as follows: 

\medskip 

$\operatorname{part}(3^t+8) = \operatorname{part}(3^t) – 
\newline 
[(3 + (3^t-3)) + (9 + (3^t-9)) + (11 + (3^t-11)) + (13 + (3^t-13))] + 
\newline 
(9) + (3 + 11 + 13) + 
\newline 
((3^t-3) + (3^t+1) + (3^t+2)) + ((3^t-9) + (3^t+3) + (3^t+6)) + 
\newline 
((3^t-11) + (3^t+4) + (3^t+7))  + ((3^t-13) + (3^t+5) + (3^t+8)).$

\medskip 

$\operatorname{part}(3^t+11) = \operatorname{part}(3^t-3) – 
\newline 
[(3 + (3^t-3)) + (9 + (3^t-9)) + (12 + (3^t-12)) + 
(13 + (3^t-13)) + (14 + (3^t-14)) + (15 + (3^t-15))] + 
\newline  
(12 + 15) + (13 + 14) + (3) + (9) + 
\newline 
((3^t-3) + 3^t + (3^t+3)) + ((3^t-15) + (3^t+6) + (3^t+9)) + 
((3^t-14) + (3^t+4) + (3^t+10)) + ((3^t-13) + (3^t+2) + (3^t+11)) + 
((3^t-12) + (3^t+5) + (3^t+7)) + ((3^t-9) + (3^t+1) + (3^t+8)).$

\medskip

Recall that there exists a unique 
3-good partition for $[3^t]$ and exactly two for $[3^t-3]$,   
and note that these partitions contain 
all sets that were subtracted in the above two examples. 

\medskip

This reduction allows us to extend 
Theorem \ref{t-main} (iii) as follows \cite{Pol25}:  

\medskip

$n = 3^t + \frac{3^\ell + 3^j}{2} - 1$, where  
$0 < j \leq \ell < t.$

\medskip  

It is easily seen that case $j=1$ corresponds to Theorem \ref{t-main} (iii) and 
that  $n \equiv 2 \pmod 3$  for all feasible $\ell, j,$ and $t$. 
The first numbers not covered by the above statements are 

\smallskip 

$n = 101, 104, 110, 113, 116, 119; 263, 266, 275, 278, 281, 290, 293, 299, \dots$ 

\medskip

It remains open 
if the generalized reduction is sufficient to prove Conjecture \ref{cnj-main}.

\section*{Acknowledgments:} 
This paper is an output of the research project (HSE-BR-2025-024) 
implemented as part of the Basic Research Program at HSE University.
We thank Oleg Polubasoff for many important suggestions and valuable insights, some of which we include in the last two sections with his kind permission. 
We are also thankful to Dmitry  Rybin for his interest in our work and important computations, 
verifying our conjectures up to $n=844$.

\end{document}